\documentclass{amsart}
\usepackage{delimset, amssymb, mathtools}
\usepackage{enumitem}
\setlist{itemsep=4pt, topsep=0pt, leftmargin=17pt, listparindent=11pt}

\usepackage{xcolor}
\definecolor{BIT}{cmyk}{1, 0, 1, 0}
\definecolor{PKU}{cmyk}{0, 1, 1, .45}
\usepackage[colorlinks, citecolor=BIT, linkcolor=PKU, pagebackref, ocgcolorlinks]{hyperref}

\newtheorem{theorem}{Theorem}[section]
\newtheorem{conjecture}[theorem]{Conjecture}
\newtheorem{corollary}[theorem]{Corollary}
\newtheorem{lemma}[theorem]{Lemma}
\newtheorem{proposition}[theorem]{Proposition}
\numberwithin{equation}{section}

\usepackage[capitalize]{cleveref}
\crefname{def}{Def.}{Defs.}
\Crefname{def}{Definition}{Definitions}
\creflabelformat{def}{~\upshape(#2#1#3)}

\allowbreak
\allowdisplaybreaks

\usepackage[numbers, sort&compress, nonamebreak, merge, elide, longnamesfirst]{natbib}

\usepackage{tikz}
\usetikzlibrary{decorations.pathreplacing, calc, positioning, arrows.meta}
\tikzset{
vertex/.style={shape=circle, minimum size=0.9mm, ball color=black, inner sep=0.5},
ellipsis/.style={shape=circle, fill, inner sep=.5},
range.lr/.style={|{Stealth}-{Stealth}|, thin},
range.l/.style={-{Stealth}|, thin},
range.r/.style={|{Stealth}-, thin},
}

\title[]{All cycle-chords are $e$-positive}
\author[D.G.L.~Wang]{David G.L. Wang$^*$}
\address[David G.L. Wang]{School of Mathematics and Statistics \& MIIT Key Laboratory of Mathematical Theory and Computation in Information Security, Beijing Institute of Technology, Beijing 102400, P.\ R.\ China.}
\email{glw@bit.edu.cn}
\thanks{$^*$Wang is supported by the NSFC (Grant No.\ 12171034).}
\begin{document}

\bibliographystyle{abbrvnat}

\begin{abstract}
We establish the $e$-positivity of cycle-chord graphs by using the composition method which is developed by Zhou and the author recently. Our method is simpler than the $(e)$-positivity approach which is used for handling cycle-chords with girth at most $4$. We also provide a combinatorial interpretation of the $e$-coefficients, and conjecture that theta graphs are $e$-positive.
\end{abstract}
\maketitle

\section{Introduction}

In 1995, \citet{Sta95} introduced the concept of 
the \emph{chromatic symmetric function} of a graph $G$
as the symmetric function
\[
X_G=\sum_{\kappa\colon V(G)\to\{1,2,\dots\}}
x_{\kappa(v_1)}x_{\kappa(v_2)}\dotsm,
\]
where the sum runs over proper colorings $\kappa$ of $G$.
It is a generalization of Birkhoff's chromatic
symmetric polynomials in the study of the $4$-color problem.
A graph $G$ is said to be \emph{$e$-positive} 
if the $e$-expansion of $X_G$ 
has only nonnegative coefficients. Here the symbol $e$ stands for
the elementary symmetric function basis
of the algebra of symmetric functions.

A leading motivation in this field 
is Stanley--Stembridge's conjecture~\cite[Conjecture 5.5]{SS93}
which asserts the $e$-positivity of the incomparability graph of 
$(3+1)$-free posets,  see also \cite[Conjecture~5.1]{Sta95}.
It is equivalent to 
the $e$-positivity of unit interval graphs,
due to the reduction of \citet{Gua13X}.
In October 2024, 
\citet{Hik24X} resolved this conjecture. 
He proved that each $e_\lambda$-coefficient of a
unit interval graph equals $\lambda!$ times
the probability of constructing 
a standard Young tableau of shape $\lambda$ 
from the empty one in a Markov chain 
according to an explicit rule of his.
\citeauthor{Hik24X} claimed that 
the probabilistic interpretation comes from 
a new geometric realization of chromatic symmetric functions
found by \citet{Kato24X}.
Common unit interval graphs include $K$-chains
such like paths, lollipops and barbells,
see \citet{Dv18,GS01,Tom24,WZ24X},
and melting lollipops, see \citet{HNY20} and \cite{Tom24}. 

In his seminal paper \cite{Sta95}, 
\citeauthor{Sta95} made an initial exploration of 
the $e$-positivity of a general graph, 
in which he stated ``a complete characterization 
of such graphs appears hopeless'';
here ``such graphs'' means $e$-positive graphs.
Indeed, only a few families of graphs that are not unit interval graphs 
were proved to be $e$-positive.
These families include 
spiders of the forms $S(a,\,a-1,\,1)$ and $S(4a+2,\,2a,\,1)$, 
see \citet{DFv20,TWW24X},
and hat chains which include cycles, tadpoles, 
generalized bulls, and hats, 
see \citet{TV24X} and \cite{WZ24X}. 

To prove the $e$-positivity of a chromatic symmetric function 
often requires specific algebraic combinatorial skills.
In this paper, 
we establish the $e$-positivity of cycle-chord graphs
by using the composition method that was recently developed by
Zhou and the author~\cite{WZ24X}.
A cycle-chord is not a unit interval graph 
as long as it has at least~$5$ vertices.
To the best of our knowledge,
only cycle-chords with girth at most $4$ is known to be $e$-positive,
see \citet{WW23-JAC}. 
The proof therein used the momentous $(e)$-positivity approach
of \citet{GS01}. 

This paper is organized as follows. 
In \cref{sec:pre},
we recall known results and tools that will be of necessary use.
In \cref{sec:main}, we establish the $e$-positivity of cycle-chords
by presenting a positive $e_I$-expansion, see \cref{thm:CC}.
We also provide a combinatorial interpretation for the $e$-coefficients,
see \cref{prop:combin}.
\Cref{sec:conj} is devoted to the $e$-positivity conjecture 
for the family of theta graphs, which
contains cycle-chords and are not unit interval graphs, see Conjecture~\ref{conj:theta}.

\section{Preliminary}\label{sec:pre}

\Citet[Theorem 3.1, Corollaries 3.2 and 3.3]{OS14} 
discovered a beautiful modular relation, 
called the \emph{triple-deletion property}, 
for chromatic symmetric functions of graphs.

\begin{proposition}[\citeauthor{OS14}]\label{prop:3del}
Let $G=(V,E)$ be a graph which contains $3$ vertices with no edges among.
Denote by $e_1$, $e_2$ and $e_3$ the edges 
to be used to link these vertices together.
For any set $S\subseteq \{1,2,3\}$, 
denote by $G_S$ the graph with vertex set~$V$
and edge set $E\cup\{e_j\colon j\in S\}$.
Then 
\[
X_{G_{12}}
=X_{G_1}+X_{G_{23}}-X_{G_3}
\quad\text{and}\quad
X_{G_{123}}
=X_{G_{13}}+X_{G_{23}}-X_{G_3}.
\]
\end{proposition}

A \emph{composition} $I$ of $n$ 
is a sequence $i_1i_2\dotsm i_z$
of positive integers whose sum is $n$, denoted $I\vDash n$.
A \emph{partition} $\lambda$ of $n$ 
is a multiset of positive integers with sum $n$, 
denoted $\lambda\vdash n$. 
It is common to arrange the parts of a partition
in the decreasing order so that it is written as a sequence as well.
Following \cite{WZ24X}, we denote by $e_I$
the elementary symmetric function indexed 
by the partition that is obtained by rearranging the parts of $I$.
Denote by $P_n$ the \emph{path} graph on $n$ vertices.
\citet[Table~1]{SW16} discovered a captivating formula
for the chromatic symmetric functions of paths.
\begin{proposition}[\citeauthor{SW16}]\label{prop:path}
We have $X_{P_n}
=\sum_{I\vDash n}
w_I
e_I$
for any $n\ge 1$, where
\begin{equation}\label[def]{def:w}
w_I
=i_1(i_2-1)(i_3-1)\dotsm(i_z-1)
\quad\text{if $I=i_1i_2\dotsm i_z$}.
\end{equation}
\end{proposition}
Denote by $C_n$ the \emph{cycle} graph on $n$ vertices. 
\citet[Corollary 6.2]{Ell17} 
gave a formula for the chromatic quasisymmetric functions of cycles, 
whose
$t=1$ specialization is equally elegant.

\begin{proposition}[\citeauthor{Ell17}]\label{prop:cycle}
We have 
$
X_{C_n}
=\sum_{I\vDash n}
(i_1-1)
w_I
e_I$
for $n\ge 2$.
\end{proposition}

For $m\ge 2$ and $l\ge 0$,
the \emph{tadpole} $C_m^l$
is the graph on $m+l$ vertices obtained by 
connecting a vertex on the cycle~$C_m$
and an end of the path $P_l$,
see the left illustration in \cref{fig:tadpole:CC}.
\begin{figure}[htbp]
\begin{tikzpicture}[decoration=brace]
\draw (0,0) node {$C_m$};
\node[vertex] (c) at (0: 1) {};
\node[vertex] at (30: 1) {};
\node[vertex] at (-30: 1) {};

\node (p1) at ($ (c) + (1, 0) $) [vertex] {};
\node[ellipsis] at ($ (p1) + (1-.2, 0) $) {};
\node[ellipsis] (e2) at ($ (p1) + (1, 0) $) {};
\node[ellipsis] at ($ (p1) + (1+.2, 0) $) {};
\node (p2) at ($ (p1) + (2, 0) $) [vertex] {};
\draw [decorate] (4, -.2) -- (1, -.2);
\node[below=9pt, font=\footnotesize] at (2.5, 0) {length $l$};

\draw (-75: 1) arc (-75: 75: 1);
\draw (c) -- ($(p1) + (1-2*.2, 0)$);
\draw (p2) -- ($(p2) - (1-2*.2, 0)$);

\begin{scope}[xshift=8cm]
\begin{scope}[rotate=90, scale=.4]
\coordinate (-2) at (-2.05, 0.6);
\coordinate (-1) at (-1.5, 2.1);
\coordinate (0)  at (0, 3);
\coordinate (1)  at (1.5, 2.1);
\coordinate (2)  at (2.05, 0.6);
\coordinate (3)  at (0, -3);
\coordinate (a)  at (-2.08, -.8); 
\coordinate (b)  at (-2.03, -1.05);
\coordinate (c)  at (-1.95, -1.3);
\coordinate (d)  at (2.08, -.8); 
\coordinate (e)  at (2.03, -1.05);
\coordinate (f)  at (1.95, -1.3);
\draw (0)--(3);
\draw (0) .. controls (1.5, 3) and (2.2, 1.1) .. (2.1,-.5);
\draw (0) .. controls (-1.5, 3) and (-2.2,1.1) .. (-2.1,-.5);
\draw (3) .. controls (-.7, -3) and (-1.6,-2.3) .. (-1.84,-1.55);
\draw (3) .. controls (.7, -3) and (1.6,-2.3) .. (1.84,-1.55);
\foreach \e in {-2,...,3} 
\shade[vertex](\e) circle(.12);
\draw (0,4.2) node{$v_0$};
\draw(-1.7, 3.5) node{$v_1$};
\draw(1.9, 4.2) node{$v_{a+b-1}$};
\draw (0,-4.2)node{$v_a$}; 
\foreach \e in {a, b, c, d, e, f}
\node[ellipsis] at (\e) {};
\end{scope} 
\end{scope}
\end{tikzpicture}
\caption{The tadpole $C_m^l$
and the cycle-chord $\mathrm{CC}_{ab}$.}\label{fig:tadpole:CC}
\end{figure}
We use the notation $C_m^l$ following \citet{QTW24X},
in which the authors study reductions 
of the chromatic symmetric functions 
of certain \emph{conjoined graphs}.
We say that a composition $I=i_1\dotsm i_z$
has \emph{modulus} $n$ if $i_1+\dots+i_z=n$, denoted $\abs{I}=n$.
The \emph{$a$-surplus} of $I$ is the number
\begin{equation}\label[def]{def:Theta+}
\Theta_I^+(a)
=\min\{\abs{i_1\dotsm i_k}\colon
0\le k\le z,\
\abs{i_1\dotsm i_k}\ge a\}
-a,
\quad
\text{for any $a\le n$.}
\end{equation}
Similarly, 
the \emph{$a$-deficiency} of $I$ is the number
\begin{equation}\label[def]{def:Theta-}
\Theta_I^-(a)
=a-\max\{\abs{i_1\dotsm i_k}\colon
0\le k\le z,\
\abs{i_1\dotsm i_k}\le a\},
\quad
\text{for any $a\ge 0$.}
\end{equation}
From definition, we see that $\Theta_I^\pm(a)\ge 0$,
and that these functions can be expressed in each other.

\begin{lemma}[\citeauthor{WZ24X}]\label{lem:Theta:+-}
For any composition $I\vDash n$
and any real number $0\le a\le n$,
we have 
$\Theta_I^-(a)
=\Theta_{\overline{I}}^+(n-a)$,
where the notation $\overline{I}$ stands for the reversal of $I$.
\end{lemma}

\begin{theorem}[\citeauthor{WZ24X}]\label{thm:tadpole}
For any $0\le l\le n-2$, we have
$X_{C_{n-l}^l}
=\sum_{I\vDash n}
\Theta_I^+(l+1)
w_I
e_I$.
\end{theorem}
It is direct to see that \cref{thm:tadpole}
reduces to Propositions~\ref{prop:path} and \ref{prop:cycle} respectively
when $l=n-2$ and $l=0$.

\section{Establishing the $e$-positivity of cycle-chords}\label{sec:main}

For $a, b\ge 1$,
the \emph{cycle-chord} $\mathrm{CC}_{ab}$
is the graph obtained by adding the
chord $v_0v_a$ in the cycle~$C_{a+b}$ 
on the vertices $v_0$, $v_1$, $\dotsc$, $v_{a+b-1}$; 
see the right figure in \cref{fig:tadpole:CC}.

\begin{lemma}\label{lem:CC}
Let $a,b\ge 1$ and $n=a+b$. Then
$X_{\mathrm{CC}_{ab}}
=\sum_{I\vDash n} c_I w_I e_I$,
where 
\[
c_I=\sum_{i=1}^b
\Theta_I^+(i)
-\sum_{i=1}^{b-1}
\Theta_{\overline{I}}^-(i).
\]
\end{lemma}
\begin{proof}
Let $G=\mathrm{CC}_{ab}$.
For $b=1$,
the desired formula reduces to Proposition~\ref{prop:cycle}.
Let $b\ge 2$.
We label its vertices in the way 
as the right illustration in \cref{fig:tadpole:CC}.
Taking $(e_1, e_2, e_3)=(v_0v_a,\, v_0v_{a+b-1},\,  v_av_{a+b-1})$, 
we obtain from the first identity of Proposition~\ref{prop:3del} that
\[
X_{\mathrm{CC}_{ab}}
=X_{\mathrm{CC}_{(a+1)(b-1)}}
+X_{C_{a+1}^{b-1}}
-X_{C_b^a}.
\]
Using it iteratedly, we can deduce that
\begin{align*}
X_{\mathrm{CC}_{ab}}
&=
\brk2{
X_{\mathrm{CC}_{(a+2)(b-2)}}
+X_{C_{a+2}^{b-2}}
-X_{C_{b-1}^{a+1}}
}
+X_{C_{a+1}^{b-1}}
-X_{C_b^a}
=\dotsm\\
&=X_{\mathrm{CC}_{(n-1)1}}
+\sum_{i=1}^{b-1}
X_{C_{a+i}^{b-i}}
-\sum_{i=0}^{b-2}
X_{C_{b-i}^{a+i}}.
\end{align*}
By merging $\mathrm{CC}_{(n-1)1}=C_n=C_{a+b}^{0}$ 
into the first sum, we obtain
\[
X_G
=\sum_{i=1}^b
X_{C_{a+i}^{b-i}}
-\sum_{i=0}^{b-2}
X_{C_{b-i}^{a+i}}.
\]
By \cref{thm:tadpole}, we can infer that
\begin{align*}
X_G
&=\sum_{i=1}^b
\sum_{I\vDash n}
\Theta_I^+(b-i+1)
w_I
e_I
-\sum_{i=0}^{b-2}
\sum_{I\vDash n}
\Theta_I^+(a+i+1)
w_I
e_I
=
\sum_{I\vDash n}
c_I
w_I
e_I,
\end{align*}
where
\[
c_I
=
\sum_{i=1}^b
\Theta_I^+(i)
-\sum_{i=a+1}^{a+b-1}
\Theta_I^+(i)
=\sum_{i=1}^b
\Theta_I^+(i)
-\sum_{i=1}^{b-1}
\Theta_{\overline{I}}^-(i),
\]
in virtue of Lemma~\ref{lem:Theta:+-}.
\end{proof}

Now we make some preparations for stating \cref{thm:CC}.
Let $I=i_1\dotsm i_z\vDash n$. 
Define $p$ and $s$ by the conditional equation
\begin{equation}\label[def]{def:ps}
b
=\abs{i_1\dotsm i_{p-1}}+s,
\quad\text{with $1\le p\le z$ and $1\le s\le i_p$}.
\end{equation}
They are well defined since $b<n$.
Define $q$ and $t$ by the conditional equation
\begin{equation}\label[def]{def:qt}
b
=\abs{i_2\dotsm i_q}+t,
\quad\text{with $1\le q\le z$, $1\le t\le i_{q+1}$ and $i_{z+1}=i_1$}.
\end{equation}
They are well defined since when $q=z$,
\cref{def:qt} implies
$t=b-(n-i_1)=i_1-a<i_1$.
We list some properties 
that the numbers $p$, $q$, $s$ and $t$ satisfy.

%: lem:psqt
\begin{lemma}\label{lem:psqt}
Let $I\vDash n$. Then we have the following.
\begin{enumerate}
\item\label[itm]{itm:q>=p-1}
$q\ge p-1$, and
$\abs{i_pi_{p+1}\dotsm i_q}+t=i_1+s$.
\item\label[itm]{itm:q=p-1:iff}
$q=p-1
\iff
t=i_1+s
\iff
i_{q+1}-t=i_p-s-i_1
\iff 
i_1\le i_p-s$.
\item\label[itm]{itm:ip-s}
$i_p-s=\Theta_{\overline{I}}^-(a)$.
\item\label[itm]{itm:q=z:iff}
$a-i_1=\abs{i_{q+1}\dotsm i_z}-t$.
As a consequence,
$q=z
\iff
t=i_1-a
\iff i_1>a$.
\end{enumerate}
\end{lemma}
\begin{proof} 
We proceed one by one.

\noindent\eqref{itm:q>=p-1}
Assume to the contrary that $q\le p-2$. Then $p\ge q+2\ge 3$.
By \cref{def:ps,def:qt},
\[
\abs{i_1\dotsm i_{p-1}}+s
=b+1
=\abs{i_2\dotsm i_q}+t
\le
\abs{i_2\dotsm i_{q+1}}
\le
\abs{i_2\dotsm i_{p-1}},
\]
which is absurd.
This proves $q\ge p-1$.
The difference of \cref{def:qt,def:ps} 
yields the identity.

\noindent\eqref{itm:q=p-1:iff}
First, if $q=p-1$, then $t=i_1+s$ by \eqref{itm:q>=p-1}.
Second, if $t=i_1+s$, then $q+1=p$ by \eqref{itm:q>=p-1}, 
and the third identity holds.
Thirdly,  
if the third identity holds, 
then $i_1\le i_p-s$ since $t\le i_{q+1}$ by premise.
At last, if $i_1\le i_p-s$,
then $q=p-1$,
since otherwise \eqref{itm:q>=p-1} would 
imply that $q\ge p$ and 
$t=i_1-(i_p-s)-\abs{i_{p+1}\dotsm i_q}\le 0$,
which is absurd.

\noindent\eqref{itm:ip-s}
Clear from \cref{def:ps}.

\noindent\eqref{itm:q=z:iff}
By \cref{def:qt},
\begin{equation}\label{pf:qt}
a-i_1
=n-\abs{i_2\dotsm i_q}-t-i_1
=\abs{i_{q+1}\dotsm i_z}-t.
\end{equation}
First, if $q=z$, then \cref{pf:qt} reduces to $t=i_1-a$ as desired.
Second, if $t=i_1-a$, then $i_1>a$ since $t\ge 1$.
At last, if $i_1>a$, 
then \cref{pf:qt} implies $\abs{i_{q+1}\dotsm i_z}<t\le i_{q+1}$.
Thus $z\le q\le z$ and $q=z$.
\end{proof}

In the remaining part of this paper, 
we need \eqref{itm:q>=p-1} and \eqref{itm:q=p-1:iff} only.
We exhibit \eqref{itm:ip-s} since 
it gives a combinatorial interpretation of the difference $i_p-s$, 
and put \eqref{itm:q=z:iff} on display
since it gives a characterization for when $q=z$.
Here comes our main result.

\begin{theorem}\label{thm:CC}
For $a,b\ge 2$ and $n=a+b$, 
\[
X_{\mathrm{CC}_{ab}}
=\sum_{I=i_1i_2\dotsm i_z\vDash n}
\Delta_I(b) w_I e_I,
\]
where $w_I=i_1(i_2-1)(i_3-1)\dotsm(i_z-1)$, and
\begin{equation}\label[def]{def:DeltaIb}
\Delta_I(b)
=\begin{dcases*}
s(i_p-s-i_1),
& if $i_1\le i_p-s$, i.e., if $q=p-1$, \\
e_2(i_p-s,\,i_{p+1},\,\dots,\,i_q,\,t),
& if $i_1>i_p-s$, i.e., if $q\ge p$,
\end{dcases*}
\end{equation}
in which $e_2(x_1,\dots,x_m)=\sum_{1\le i<j\le m}x_i x_j$,
and the symbols $p,q,s,t$ are defined by 
\[
b=\abs{i_1\dotsm i_{p-1}}+s
=\abs{i_2\dotsm i_q}+t,
\]
such that $1\le p,q\le z$, $1\le s\le i_p$, and $1\le t\le i_{q+1}$,
with the convention $i_{z+1}=i_1$.
\end{theorem}
\begin{proof}
Let $G=\mathrm{CC}_{ab}$ and $n=a+b$.
We define
a transformation $\varphi$ on the set of compositions of $n$ by
\[
\varphi(i_1\dotsm i_z)
=i_1 i_z i_{z-1}\dotsm i_2.
\]
Then $\varphi$ is an involution, $w_{\varphi(I)}=w_I$
and $e_{\varphi(I)}=e_I$. 
By Lemma~\ref{lem:CC}, 
\begin{align*}
c_I+c_{\varphi(I)}
&=\brk4{\sum_{i=1}^b
\Theta_I^+(i)
-\sum_{i=1}^{b-1}
\Theta_{\overline{I}}^-(i)}
+\brk4{\sum_{i=1}^b
\Theta_{\varphi(I)}^+(i)
-\sum_{i=1}^{b-1}
\Theta_{\overline{\varphi(I)}}^-(i)}.
\end{align*}
Regrouping the four sums above, we obtain
$c_I+c_{\varphi(I)}
=\Delta_I'(b)
+\Delta_{\varphi(I)}'(b)$,
where 
\begin{equation}\label{pf:CC:DeltaIb}
\Delta_I'(b)
=\sum_{i=1}^b
\Theta_I^+(i)
-\sum_{i=1}^{b-1}
\Theta_{\overline{\varphi(I)}}^-(i).
\end{equation}
Since $\varphi$ is an involution 
on the set of compositions of $n$, we can deduce that
\begin{align*}
X_G
=\frac{1}{2}
\sum_{I\vDash n}
\brk1{
c_I+c_{\varphi(I)}}
w_I
e_I
=\frac{1}{2}
\sum_{I\vDash n}
\brk1{
\Delta_I'(b)
+\Delta_{\varphi(I)}'(b)}
w_I
e_I
=\sum_{I\vDash n}
\Delta_I'(b)
w_I
e_I.
\end{align*}
We shall show that $\Delta_I'(b)=\Delta_I(b)$.

First,
we calculate the two sums in \cref{pf:CC:DeltaIb}.
By \cref{def:Theta+,def:ps},
\begin{align}
\notag
\sum_{i=1}^b
\Theta_I^+(i)
&=
\brk1{
(i_1-1)+(i_1-2)+\dots+0}
+
\dots
+
\brk1{
(i_{p-1}-1)+(i_{p-1}-2)+\dots+0}
\\
\notag
&\qquad
+
\brk1{
(i_p-1)+(i_p-2)+\dots+(i_p-s)}
\\
\label{pf:sum:Theta+}
&=\binom{i_1}{2}
+\dots
+\binom{i_{p-1}}{2}
+si_p
-\binom{s+1}{2}.
\end{align}
It is elementary to show that
the $e_2$-function satisfies the identity
\begin{equation}\label{id}
\sum_{j=1}^k\binom{x_j}{2}
=\binom{x_1+\dots+x_k}{2}
-e_2(x_1,\dots,x_k).
\end{equation}
If we set $(x_1, \dots, x_k)=(i_1,\, \dots,\, i_{p-1},\, s)$, 
then \cref{id} becomes
\[
\sum_{j=1}^{p-1}
\binom{i_j}{2}
+\binom{s}{2}
=\binom{\abs{i_1\dotsm i_{p-1}}+s}{2}
-e_2(i_1,\,\dots,\,i_{p-1},\,s)
=
\binom{b}{2}
-e_2(i_1,\,\dots,\,i_{p-1})
-s(b-s).
\]
Substituting it into \cref{pf:sum:Theta+}
and by algebraic simplification, 
we obtain
\begin{equation}\label{pf:CC:Delta:sum1}
\sum_{i=1}^b
\Theta_I^+(i)
=
\binom{b}{2}
-e_2(i_1,\,\dots,\,i_{p-1})
-s(b-i_p).
\end{equation}
On the other hand, 
recall that $\overline{\varphi(I)}=i_2i_3\dotsm i_z i_1$.
By \cref{def:Theta-,def:qt},
and by setting
in \cref{id}
$(x_1, \dots, x_k)=(i_2,\,\dots,\,i_q,\,t)$, 
we can deduce that 
\begin{multline}\label{pf:CC:Delta:sum2}
\sum_{i=1}^{b-1}
\Theta_{\overline{\varphi(I)}}^-(i)
=
\brk1{
1+\dots+(i_2-1)}
+
\dots
+
\brk1{
1+\dots+(i_q-1)}
+
\brk1{
1+\dots+(t-1)}
\\
=
\binom{i_2}{2}
+\dots
+\binom{i_q}{2}
+\binom{t}{2}
=
\binom{b}{2}
-e_2(i_2,\,\dots,\,i_q,\,t).
\end{multline}
Substituting \cref{pf:CC:Delta:sum1,pf:CC:Delta:sum2} into 
\cref{pf:CC:DeltaIb}, we obtain
\[
\Delta_I'(b)
=e_2(i_2,\,\dots,\,i_q,\,t)
-e_2(i_1,\,\dots,\,i_{p-1})
-s(b-i_p).
\]
It remains to simplify it to \cref{def:DeltaIb}.
Since $q\ge p-1$ by Lemma~\ref{lem:psqt}, we can expand 
\begin{align*}
e_2(i_2,\,\dots,\,i_q,\,t)
&=
e_2(i_2,\,\dots,\,i_{p-1})
+
e_2(i_p,\,\dots,\,i_q,\,t)
+
\abs{i_2\dotsm i_{p-1}}
\brk1{\abs{i_p\dotsm i_q}+t}, 
\quad\text{and}\\
e_2(i_1,\,\dots,\,i_{p-1})
&=
e_2(i_2,\,\dots,\,i_{p-1})
+
i_1\abs{i_2\dotsm i_{p-1}}.
\end{align*}
Substituting them into the previous expression,  we obtain
\[
\Delta_I'(b)
=
e_2(i_p,\,\dots,\,i_q,\,t)
+
\abs{i_2\dotsm i_{p-1}}
\brk1{\abs{i_pi_{p+1}\dotsm i_q}+t-i_1}
-s(b-i_p).
\]
Then by Lemma~\ref{lem:psqt} \eqref{itm:q>=p-1} and \cref{def:ps}, 
we can deduce that
\begin{equation}\label{eq:DeltaI'b}
\Delta_I'(b)
-e_2(i_p,\,\dots,\,i_q,\,t)
=
\abs{i_2\dotsm i_{p-1}}s
-s\brk1{
\abs{i_1\dotsm i_{p-1}}+s-i_p}
=s(i_p-s-i_1).
\end{equation}
If $q=p-1$, i.e., if $i_1\le i_p-s$
by Lemma~\ref{lem:psqt} \eqref{itm:q=p-1:iff}, 
then \cref{eq:DeltaI'b} 
implies $\Delta_I'(b)=s(i_p-s-i_1)=\Delta_I(b)$.
Otherwise $q\ge p$. 
By Lemma~\ref{lem:psqt} \eqref{itm:q>=p-1}, 
\cref{eq:DeltaI'b} implies
\begin{align*}
\Delta_I'(b)
&=e_2(i_p,\,\dots,\,i_q,\,t)
-s(s+i_1-i_p)\\
&=e_2(i_{p},\,\dots,\,i_q,\,t)
-s\brk1{\abs{i_p\dotsm i_q}+t-i_p}\\
&=e_2(i_p-s,\,i_{p+1},\,\dots,\,i_q,\,t)
=\Delta_I(b).
\end{align*}
This completes the proof.
\end{proof}

We remark that
the condition 
$i_1>i_p-s$ in \cref{def:DeltaIb} cannot be replaced with $i_1\ge i_p-s$,
since the $e_2$-function is not well defined when $i_1=i_p-s$.
For any partition $\lambda\vdash n$ and any symmetric function~$X$,
we use the notation $[e_\lambda]X$ to denote the coefficient
of $e_\lambda$ in the $e$-expansion of $X$.

\begin{corollary}\label{cor:e-coeff:CC}
For any $a,b\ge 2$, the cycle-chord $\mathrm{CC}_{ab}$
is $e$-positive. Moreover,
for $n=a+b$, we have
$[e_n]X_{\mathrm{CC}_{ab}}=abn$
and
$[e_{(n-1)1}]X_{\mathrm{CC}_{ab}}=(a-1)(b-1)(n-2)$.
\end{corollary}
\begin{proof}
Let $G=\mathrm{CC}_{ab}$.
The $e$-positivity is clear since $\Delta_I(b)\ge 0$ 
by \cref{def:DeltaIb}.
Taking $I=n$ in \cref{thm:CC}, 
we obtain $z=p=q=1$ and $b=s=t$.
Hence $[e_n]X_G=e_2(n-s,\,t)w_n=abn$.
On the other hand, 
since $w_{(n-1)1}=0$, 
we set $I=1(n-1)$ in \cref{thm:CC}.
Then $z=p=q=2$ and $s=b-1$. It follows that 
$[e_{(n-1)1}]X_G=s(i_p-s-i_1)w_I=(b-1)(a-1)(n-2)$.
\end{proof}

Here is a combinatorial interpretation for the numbers $\Delta_I(b)$.

\begin{theorem}\label{prop:combin}
Let $I=i_1\dotsm i_z\vDash n$.
In the unique dissection
of the interval $(0,\,n+i_1]$ into open-closed intervals
$\mathcal I_j$ of lengths $i_1$, $i_2$, $\dots$, $i_z$, $i_1$,
we call each $\mathcal I_j$ a \emph{segment}.
Let $\mathcal B=(b,\,b+i_1]$.
\begin{enumerate}
\item
If $\mathcal B$ is contained in a single segment,
then the number $\Delta_I(b)$
equals the product of lengths 
of the two sub-segments obtained by removing
$\mathcal B$ from that segment.
\item
Otherwise, the number $\Delta_I(b)$
is evaluated by the $e_2$-function in the lengths 
of nonempty intersections $\mathcal I\cap\mathcal B$,
where $\mathcal I$ runs over the $z+1$ segments.
\end{enumerate}
As a consequence, $\Delta_I(b)=0$ if and only if 
$\mathcal B$ is contained in a single segment 
of the form $(b,y]$ or $(x,\,b+i_1]$,
where $x$ and $y$ are integers.
\end{theorem}
\begin{proof}
If $i_1\le i_p-s$,
then $\mathcal B\subseteq\mathcal I_p$,
and the sub-segments in
$\mathcal I_p\backslash \mathcal B$
have lengths $s$ and $i_p-s-i_1$, respectively, 
see \cref{fig:i1<}.
\begin{figure}[htbp]
\begin{tikzpicture}
\draw (1, 0) -- (10, 0);
\draw[range.lr] (3, .12) -- node[above] {$i_p$} ++(5,0);
\draw[range.lr] (3, -.12) -- node[below=2pt] {$s$} ++(1,0);
\draw[range.lr] (4, -.12) -- node[below] {$i_1$} ++(1,0);
\draw[range.lr] (5, -.12) -- node[below] {$i_p-s-i_1$} ++(3,0);
\draw[range.l] (1, -.8) -- node[below] {$b$} ++(3,0);
\draw[range.r] (4, -.8) -- node[below] {$a$} ++(6,0);
\draw[ultra thick, red] 
(3,0) -- (4,0) 
(5,0) -- (8,0);
\node[vertex] at (3,0) {};
\node[vertex] at (4,0) {};
\node[vertex] at (5,0) {};
\node[vertex] at (8,0) {};
\end{tikzpicture}
\caption{When $i_1\le i_p-s$, we have $\Delta_I(b)=s(i_p-s-i_1)$.}
\label{fig:i1<}
\end{figure}
If $i_1>i_p-s$,
then $q\ge p$
by Lemma~\ref{lem:psqt} \eqref{itm:q>=p-1} and \eqref{itm:q=p-1:iff}.
Moreover, the variables in the $e_2$-function
have sum $(i_p-s)+\abs{i_{p+1}\dotsm i_q}+t=i_1$,
see \cref{fig:i1>}.
\begin{figure}[htbp]
\begin{tikzpicture}
\draw (1, 0) -- (14, 0);
\draw[range.lr] (3, .12) -- node[above] {$i_p$} ++(3,0);
\draw[range.lr] (10, .12) -- node[above] {$i_{q+1}$} ++(2,0);
\draw[range.lr] (6, -.12) -- node[below] {$i_{p+1}$} ++(1,0);
\draw (8,0) node[below=2pt] {$\dotsm$};
\draw[range.lr] (9, -.12) -- node[below] {$i_q$} ++(1,0);
\draw[range.lr] (3, -.12) -- node[below=2pt] {$s$} ++(1,0);
\draw[range.lr] (4, -.12) -- node[below] {$i_p-s$} ++(2,0);
\draw[range.lr] (10, -.12) -- node[below=1pt] {$t$} ++(1,0);
\draw[range.lr] (4, -.8) -- node[below] {$i_1$} ++(7,0);
\draw[range.l] (1, -.8) -- node[below]{$b$} ++(3, 0);
\draw[ultra thick, red] (4,0) -- (11,0);
\node[vertex] at (3,0) {};
\node[vertex] at (4,0) {};
\node[vertex] at (6,0) {};
\node[vertex] at (7,0) {};
\node[vertex] at (9,0) {};
\node[vertex] at (10,0) {};
\node[vertex] at (11,0) {};
\node[vertex] at (12,0) {};
\end{tikzpicture}
\caption{When $i_1>i_p-s$, 
we have $\Delta_I(b)=e_2(i_p-s,\,i_{p+1},\,\dots,\,i_q,\,t)$.}
\label{fig:i1>}
\end{figure}
In this case, it is also possible that~$\mathcal B$ is contained 
in a single segment,
which happens if and only if $i_p=s$ and $q=p$.
In conclusion, the interpretation follows from \cref{thm:CC},
and the consequence is clear from the interpretation.
\end{proof}

In the language of general topology, 
the variables in the $e_2$-function are the lengths of members 
in the cover 
$\bigcup_{j=1}^{z+1} (\mathcal I_j\cap\mathcal B)$
of $\mathcal B$ that is induced by 
the segment cover $\bigcup_{j=1}^{z+1}\mathcal I_j$
of the interval $(0,\,n+i_1]$.
For example, consider the graph $\mathrm{CC}_{33}$
that is obtained by identifying an edge of two squares $C_4$. 
It is routine to calculate that 
\[
X_{\mathrm{CC}_{33}}
=\sum_{I\vDash 6}\Delta_I(3)w_Ie_I
=54e_6+16e_{51}+26e_{42}+2e_{222}.
\]
The $e_{42}$-coefficient~$26$ can be verified 
by \cref{prop:combin} as follows.
\begin{enumerate}
\item
For $I=24$, we have $\mathcal B=(3,5]$,
which is contained in the single segment $(2,6]$.
Since the sub-segments obtained by removing~$\mathcal B$ 
from $(2,6]$
are $(2,3]$ and $(5,6]$,
we find $\Delta_I(3)=1$.
\item
For $I=42$, we have $\mathcal B=(3,7]$, 
and the segments $\mathcal I$ are $(0,4]$, $(4,6]$ and $(6,10]$.
The nonempty intersections $\mathcal I\cap\mathcal B$
are $(3,4]$, $(4,6]$ and $(6,7]$.
Thus $\Delta_I(3)=e_2(1,2,1)=5$.
\end{enumerate}
Therefore, $[e_{42}]X_{\mathrm{CC}_{33}}=w_{24}+5w_{42}=26$.

\section{The $e$-positivity conjecture for theta graphs}\label{sec:conj}

For any partition $\lambda=\lambda_1\dotsm\lambda_l$,
let $G_\lambda$ be the graph of order $n=\abs{\lambda}-\ell(\lambda)+2$
obtained from two distinct vertices 
by linking them with pairwise disjoint paths 
of lengths $\lambda_1,\,\dots,\,\lambda_l$, respectively.
When $l=1$, $G_\lambda$ is a path;
when $l=2$, $G_\lambda$ is a cycle.
The $e$-positivity of paths and cycles are classical results due to 
\citet[Propositions~5.3 and 5.4]{Sta95}.
When $l\ge 3$,
we can suppose that $\lambda_{l-1}\ge 2$
since 
$X_{G_{\lambda_1\dotsm\lambda_{l-2}11}}
=X_{G_{\lambda_1\dotsm\lambda_{l-2}1}}$.
When $l=3$, $G_\lambda$ is called a \emph{theta} graph
and denoted $\theta_\lambda$, c.f.~\citet{Bon72}.
We pose the following $e$-positivity conjecture for theta graphs.

\begin{conjecture}\label{conj:theta}
For any integers $a\ge b\ge c\ge 1$,
the theta graph $\theta_{abc}$ is $e$-positive.
\end{conjecture}

The case $c=1$ is confirmed by \cref{thm:CC}.
The case $c=2$ was recently confirmed by \citet{CHW24X}.
In the spirit of positivity classification, 
any graph uniquely belongs to one of the three classes:
(i) $e$-positive graphs, 
(ii) Schur positive graphs that are not $e$-positive,
and (iii) graphs that are not Schur positive, see \citet{WW20}.
When $l\ge 4$, we know that each of these three families 
contains a graph of the form $G_\lambda$. 
In fact, it is routine to check that $G_{3221}$ and $G_{3321}$ 
belong to the first two classes respectively.
\citet[Proposition~1.5]{Sta98} 
showed that any Schur positive graph~$G$ is nice, 
that is, if $G$ contains a stable partition of type $\lambda$,
then $G$ contains a stable partition 
of any type $\mu$ that is dominated by $\lambda$.
As an application, the graph $G_{2^k1}$ for any $k\ge 3$ 
is not Schur positive since it is not nice;
it contains a stable partition of type $k1^2$ 
but no stable partitions of type $(k-1)21$. 
Similarly, the graph $G_{2^k}$ for any $k\ge 4$ is not nice, 
since it contains a stable partition of type $k2$ 
but no stable partitions of type $(k-1)3$.

\section*{Acknowledgement}
The author is appreciative of the hospitality 
of the Research Center for Mathematics and Interdisciplinary Sciences 
of Shandong University,
where he was visiting Professor Shuxiao Li and Professor Zhicong Lin.
He is also grateful to the encouragement of Professor Jean-Yves Thibon
when the author is visiting the LIGM 
at Université Gustave Eiffel in 2024. 
He thanks the anonymous reviewer 
for suggesting fleshing out the combinatorial interpretation.

\section*{Declarations}
All data generated or analyzed during 
this study are included in this paper.
The author has no relevant financial 
or non-financial interests to disclose.
He has no conflicts of interest to declare 
that are relevant to the content of this paper.

\bibliography{csf-CC}

\begin{thebibliography}{22}
\providecommand{\natexlab}[1]{#1}
\providecommand{\url}[1]{\texttt{#1}}
\expandafter\ifx\csname urlstyle\endcsname\relax
  \providecommand{\doi}[1]{doi: #1}\else
  \providecommand{\doi}{doi: \begingroup \urlstyle{rm}\Url}\fi

\bibitem[Bondy(1972)]{Bon72}
J.~A. Bondy.
\newblock The ``graph theory'' of the greek alphabet.
\newblock In Y.~Alavi, D.~R. Lick, and A.~T. White, editors, \emph{Graph Theory
  and Applications}, volume 303 of \emph{Lecture Notes in Math. (LNM)}, pages
  43--54, Berlin, Heidelberg, 1972. Springer Berlin Heidelberg.

\bibitem[Chen et~al.(2024)Chen, He, and Wang]{CHW24X}
L.~Chen, Y.~He, and D.~G.~L. Wang.
\newblock Clocks are $e$-positive.
\newblock arXiv: 2410.07581, 2024.

\bibitem[Dahlberg and {van Willigenburg}(2018)]{Dv18}
S.~Dahlberg and S.~{van Willigenburg}.
\newblock Lollipop and lariat symmetric functions.
\newblock \emph{SIAM J. Discrete Math.}, 32\penalty0 (2):\penalty0 1029--1039,
  2018.

\bibitem[Dahlberg et~al.(2020)Dahlberg, Foley, and {van Willigenburg}]{DFv20}
S.~Dahlberg, A.~Foley, and S.~{van Willigenburg}.
\newblock Resolving {S}tanley's $e$-positivity of claw-contractible-free
  graphs.
\newblock \emph{J. European Math. Soc.}, 22\penalty0 (8):\penalty0 2673--2696,
  2020.

\bibitem[Ellzey(2017)]{Ell17}
B.~Ellzey.
\newblock Chromatic quasisymmetric functions of directed graphs.
\newblock \emph{S\'{e}m. Lothar. Combin.}, 78B:\penalty0 Art. 74, 12, 2017.

\bibitem[Gebhard and Sagan(2001)]{GS01}
D.~D. Gebhard and B.~E. Sagan.
\newblock A chromatic symmetric function in noncommuting variables.
\newblock \emph{J. Alg. Combin.}, 13\penalty0 (3):\penalty0 227--255, 2001.

\bibitem[Guay-Paquet(2013)]{Gua13X}
M.~Guay-Paquet.
\newblock A modular relation for the chromatic symmetric functions of
  (3+1)-free posets.
\newblock arXiv: 1306.2400, 2013.

\bibitem[Hikita(2024)]{Hik24X}
T.~Hikita.
\newblock A proof of the {Stanley--Stembridge} conjecture.
\newblock arXiv: 2410.12758, 2024.

\bibitem[Huh et~al.(2020)Huh, Nam, and Yoo]{HNY20}
J.~Huh, S.-Y. Nam, and M.~Yoo.
\newblock Melting lollipop chromatic quasisymmetric functions and {S}chur
  expansion of unicellular {LLT} polynomials.
\newblock \emph{Discrete Math.}, 343\penalty0 (3):\penalty0 111728, 2020.

\bibitem[Kato(2024)]{Kato24X}
S.~Kato.
\newblock A geometric realization of the chromatic symmetric function of a unit
  interval graph.
\newblock arXiv: 2410.12231, 2024.

\bibitem[Orellana and Scott(2014)]{OS14}
R.~Orellana and G.~Scott.
\newblock Graphs with equal chromatic symmetric functions.
\newblock \emph{Discrete Math.}, 320:\penalty0 1--14, 2014.

\bibitem[Qi et~al.(2024)Qi, Tang, and Wang]{QTW24X}
E.~Y.~J. Qi, D.~Q.~B. Tang, and D.~G.~L. Wang.
\newblock Chromatic symmetric functions of conjoined graphs.
\newblock arXiv: 2406.01418, 2024.

\bibitem[Shareshian and Wachs(2016)]{SW16}
J.~Shareshian and M.~L. Wachs.
\newblock Chromatic quasisymmetric functions.
\newblock \emph{Adv. Math.}, 295:\penalty0 497--551, 2016.

\bibitem[Stanley(1995)]{Sta95}
R.~P. Stanley.
\newblock A symmetric function generalization of the chromatic polynomial of a
  graph.
\newblock \emph{Adv. Math.}, 111\penalty0 (1):\penalty0 166--194, 1995.

\bibitem[Stanley(1998)]{Sta98}
R.~P. Stanley.
\newblock Graph colorings and related symmetric functions: ideas and
  applications: a description of results, interesting applications, \& notable
  open problems.
\newblock \emph{Discrete Math.}, 193\penalty0 (1-3):\penalty0 267--286, 1998.

\bibitem[Stanley and Stembridge(1993)]{SS93}
R.~P. Stanley and J.~R. Stembridge.
\newblock On immanants of {J}acobi-{T}rudi matrices and permutations with
  restricted position.
\newblock \emph{J. Combin. Theory Ser. A}, 62\penalty0 (2):\penalty0 261--279,
  1993.

\bibitem[Tang et~al.(2024)Tang, Wang, and Wang]{TWW24X}
D.~Q.~B. Tang, D.~G.~L. Wang, and M.~M.~Y. Wang.
\newblock The spiders ${S}(4m+2,\,2m,\,1)$ are $e$-positive.
\newblock arXiv: 2405.04915, 2024.

\bibitem[Tom(2024)]{Tom24}
F.~Tom.
\newblock A signed $e$-expansion of the chromatic symmetric function and some
  new $e$-positive graphs.
\newblock In \emph{Proceedings of the 36th Conference on Formal Power Series
  and Algebraic Combinatorics, S\'{e}m. Lothar. Combin.}, volume 91B, Article
  \#48, 12 pp., Bochum, 2024.

\bibitem[Tom and Vailaya(2024)]{TV24X}
F.~Tom and A.~Vailaya.
\newblock Adjacent cycle-chains are $e$-positive.
\newblock arXiv: 2410.21762, 2024.

\bibitem[Wang and Wang(2020)]{WW20}
D.~G.~L. Wang and M.~M.~Y. Wang.
\newblock A combinatorial formula for the {S}chur coefficients of chromatic
  symmetric functions.
\newblock \emph{Discrete Appl. Math.}, 285:\penalty0 621--630, 2020.

\bibitem[Wang and Wang(2023)]{WW23-JAC}
D.~G.~L. Wang and M.~M.~Y. Wang.
\newblock The $e$-positivity of two classes of cycle-chord graphs.
\newblock \emph{J. Alg. Combin.}, 57\penalty0 (2):\penalty0 495--514, 2023.

\bibitem[Wang and Zhou(2025)]{WZ24X}
D.~G.~L. Wang and J.~Z.~F. Zhou.
\newblock Composition method for chromatic symmetric functions: Neat
  noncommutative analogs.
\newblock arXiv: 2401.01027, to be appeared in Adv. Appl. Math., 2025.

\end{thebibliography}

\end{document}